\documentclass{article}
\usepackage{amssymb,amsmath,theorem,euscript}

\usepackage{upgreek}

\input{epsf}

\newcounter{sec}

\def\sm{\smallskip}

\newcounter{punct}[sec]

\def\punct{\refstepcounter{punct}{\arabic{sec}.\arabic{punct}.  }}

\def\COUNTERS{\addtocounter{sec}{1}
              \setcounter{punct}{0}
          \setcounter{equation}{0}
          \setcounter{theorem}{0}
                  }

\newtheorem{ttheorem}{Theorem}
\newtheorem{llemma}{Lemma}

\begin{document}

 \def\ov{\overline}
\def\wt{\widetilde}
\def\wh{\widehat}
 \newcommand{\sgn}{\mathop {\mathrm {sign}}\nolimits}
 \newcommand{\rk}{\mathop {\mathrm {rk}}\nolimits}
\newcommand{\Aut}{\mathop {\mathrm {Aut}}\nolimits}
\newcommand{\Out}{\mathop {\mathrm {Out}}\nolimits}
\newcommand{\Abs}{\mathop {\mathrm {Abs}}\nolimits}
\renewcommand{\Re}{\mathop {\mathrm {Re}}\nolimits}
\renewcommand{\Im}{\mathop {\mathrm {Im}}\nolimits}
 \newcommand{\tr}{\mathop {\mathrm {tr}}\nolimits}
  \newcommand{\Hom}{\mathop {\mathrm {Hom}}\nolimits}
   \newcommand{\diag}{\mathop {\mathrm {diag}}\nolimits}
   \newcommand{\supp}{\mathop {\mathrm {supp}}\nolimits}
 \newcommand{\im}{\mathop {\mathrm {im}}\nolimits}
 \newcommand{\grad}{\mathop {\mathrm {grad}}\nolimits}
  \newcommand{\sgrad}{\mathop {\mathrm {sgrad}}\nolimits}
 \newcommand{\rot}{\mathop {\mathrm {rot}}\nolimits}
  \renewcommand{\div}{\mathop {\mathrm {div}}\nolimits}

\def\Br{\mathrm {Br}}
\def\Vir{\mathrm {Vir}}

 \def\Ham{\mathrm {Ham}}
\def\SL{\mathrm {SL}}
\def\Pol{\mathrm {Pol}}
\def\SU{\mathrm {SU}}
\def\GL{\mathrm {GL}}
\def\U{\mathrm U}
\def\OO{\mathrm O}
 \def\Sp{\mathrm {Sp}}
  \def\Ad{\mathrm {Ad}}
 \def\SO{\mathrm {SO}}
\def\SOS{\mathrm {SO}^*}
 \def\Diff{\mathrm{Diff}}
 \def\Vect{\mathfrak{Vect}}
\def\PGL{\mathrm {PGL}}
\def\PU{\mathrm {PU}}
\def\PSL{\mathrm {PSL}}
\def\Symp{\mathrm{Symp}}
\def\Cont{\mathrm{Cont}}
\def\End{\mathrm{End}}
\def\Mor{\mathrm{Mor}}
\def\Aut{\mathrm{Aut}}
 \def\PB{\mathrm{PB}}
\def\Fl{\mathrm {Fl}}
\def\Symm{\mathrm {Symm}} 
 \def\Herm{\mathrm {Herm}} 
  \def\SDiff{\mathrm {SDiff}} 
 
 \def\cA{\mathcal A}
\def\cB{\mathcal B}
\def\cC{\mathcal C}
\def\cD{\mathcal D}
\def\cE{\mathcal E}
\def\cF{\mathcal F}
\def\cG{\mathcal G}
\def\cH{\mathcal H}
\def\cJ{\mathcal J}
\def\cI{\mathcal I}
\def\cK{\mathcal K}
 \def\cL{\mathcal L}
\def\cM{\mathcal M}
\def\cN{\mathcal N}
 \def\cO{\mathcal O}
\def\cP{\mathcal P}
\def\cQ{\mathcal Q}
\def\cR{\mathcal R}
\def\cS{\mathcal S}
\def\cT{\mathcal T}
\def\cU{\mathcal U}
\def\cV{\mathcal V}
 \def\cW{\mathcal W}
\def\cX{\mathcal X}
 \def\cY{\mathcal Y}
 \def\cZ{\mathcal Z}
\def\0{{\ov 0}}
 \def\1{{\ov 1}}
 \def\frA{\mathfrak A}
 \def\frB{\mathfrak B}
\def\frC{\mathfrak C}
\def\frD{\mathfrak D}
\def\frE{\mathfrak E}
\def\frF{\mathfrak F}
\def\frG{\mathfrak G}
\def\frH{\mathfrak H}
\def\frI{\mathfrak I}
 \def\frJ{\mathfrak J}
 \def\frK{\mathfrak K}
 \def\frL{\mathfrak L}
\def\frM{\mathfrak M}
 \def\frN{\mathfrak N} \def\frO{\mathfrak O} \def\frP{\mathfrak P} \def\frQ{\mathfrak Q} \def\frR{\mathfrak R}
 \def\frS{\mathfrak S} \def\frT{\mathfrak T} \def\frU{\mathfrak U} \def\frV{\mathfrak V} \def\frW{\mathfrak W}
 \def\frX{\mathfrak X} \def\frY{\mathfrak Y} \def\frZ{\mathfrak Z} \def\fra{\mathfrak a} \def\frb{\mathfrak b}
 \def\frc{\mathfrak c} \def\frd{\mathfrak d} \def\fre{\mathfrak e} \def\frf{\mathfrak f} \def\frg{\mathfrak g}
 \def\frh{\mathfrak h} \def\fri{\mathfrak i} \def\frj{\mathfrak j} \def\frk{\mathfrak k} \def\frl{\mathfrak l}
 \def\frm{\mathfrak m} \def\frn{\mathfrak n} \def\fro{\mathfrak o} \def\frp{\mathfrak p} \def\frq{\mathfrak q}
 \def\frr{\mathfrak r} \def\frs{\mathfrak s} \def\frt{\mathfrak t} \def\fru{\mathfrak u} \def\frv{\mathfrak v}
 \def\frw{\mathfrak w} \def\frx{\mathfrak x} \def\fry{\mathfrak y} \def\frz{\mathfrak z} \def\frsp{\mathfrak{sp}}
 \def\bfa{\mathbf a} \def\bfb{\mathbf b} \def\bfc{\mathbf c} \def\bfd{\mathbf d} \def\bfe{\mathbf e} \def\bff{\mathbf f}
 \def\bfg{\mathbf g} \def\bfh{\mathbf h} \def\bfi{\mathbf i} \def\bfj{\mathbf j} \def\bfk{\mathbf k} \def\bfl{\mathbf l}
 \def\bfm{\mathbf m} \def\bfn{\mathbf n} \def\bfo{\mathbf o} \def\bfp{\mathbf p} \def\bfq{\mathbf q} \def\bfr{\mathbf r}
 \def\bfs{\mathbf s} \def\bft{\mathbf t} \def\bfu{\mathbf u} \def\bfv{\mathbf v} \def\bfw{\mathbf w} \def\bfx{\mathbf x}
 \def\bfy{\mathbf y} \def\bfz{\mathbf z} \def\bfA{\mathbf A} \def\bfB{\mathbf B} \def\bfC{\mathbf C} \def\bfD{\mathbf D}
 \def\bfE{\mathbf E} \def\bfF{\mathbf F} \def\bfG{\mathbf G} \def\bfH{\mathbf H} \def\bfI{\mathbf I} \def\bfJ{\mathbf J}
 \def\bfK{\mathbf K} \def\bfL{\mathbf L} \def\bfM{\mathbf M} \def\bfN{\mathbf N} \def\bfO{\mathbf O} \def\bfP{\mathbf P}
 \def\bfQ{\mathbf Q} \def\bfR{\mathbf R} \def\bfS{\mathbf S} \def\bfT{\mathbf T} \def\bfU{\mathbf U} \def\bfV{\mathbf V}
 \def\bfW{\mathbf W} \def\bfX{\mathbf X} \def\bfY{\mathbf Y} \def\bfZ{\mathbf Z} \def\bfw{\mathbf w}
 \def\R {{\mathbb R }} \def\C {{\mathbb C }} \def\Z{{\mathbb Z}} \def\H{{\mathbb H}} \def\K{{\mathbb K}}
 \def\N{{\mathbb N}} \def\Q{{\mathbb Q}} \def\A{{\mathbb A}} \def\T{\mathbb T} \def\P{\mathbb P} \def\G{\mathbb G}
 \def\bbA{\mathbb A} \def\bbB{\mathbb B} \def\bbD{\mathbb D} \def\bbE{\mathbb E} \def\bbF{\mathbb F} \def\bbG{\mathbb G}
 \def\bbI{\mathbb I} \def\bbJ{\mathbb J} \def\bbL{\mathbb L} \def\bbM{\mathbb M} \def\bbN{\mathbb N} \def\bbO{\mathbb O}
 \def\bbP{\mathbb P} \def\bbQ{\mathbb Q} \def\bbS{\mathbb S} \def\bbT{\mathbb T} \def\bbU{\mathbb U} \def\bbV{\mathbb V}
 \def\bbW{\mathbb W} \def\bbX{\mathbb X} \def\bbY{\mathbb Y} \def\kappa{\varkappa} \def\epsilon{\varepsilon}
 \def\phi{\varphi} \def\le{\leqslant} \def\ge{\geqslant}

\def\UU{\bbU}
\def\Mat{\mathrm{Mat}}
\def\tto{\rightrightarrows}

\def\F{\mathbf{F}}

\def\Gms{\mathrm {Gms}}
\def\Ams{\mathrm {Ams}}
\def\Isom{\mathrm {Isom}}

\def\Gr{\mathrm{Gr}}

\def\graph{\mathrm{graph}}

\def\la{\langle}
\def\ra{\rangle}

\def\lla{\la\!\la}
\def\rra{\ra\!\ra}


 \def\ov{\overline}
\def\wt{\widetilde}

\renewcommand{\Re}{\mathop {\mathrm {Re}}\nolimits}
\def\Br{\mathrm {Br}}

  \def\Match{\mathrm {Match}}
 \def\Isom{\mathrm {Isom}}
 \def\Hier{\mathrm {Hier}}
\def\SL{\mathrm {SL}}
\def\SU{\mathrm {SU}}
\def\GL{\mathrm {GL}}
\def\U{\mathrm U}
\def\OO{\mathrm O}
 \def\Sp{\mathrm {Sp}}
  \def\GLO{\mathrm {GLO}}
 \def\SO{\mathrm {SO}}
\def\SOS{\mathrm {SO}^*}
 \def\Diff{\mathrm{Diff}}
 \def\Vect{\mathfrak{Vect}}
\def\PGL{\mathrm {PGL}}
\def\PU{\mathrm {PU}}
\def\PSL{\mathrm {PSL}}
\def\Symp{\mathrm{Symp}}
\def\ASymm{\mathrm{Asymm}}
\def\Asymm{\mathrm{Asymm}}
\def\Gal{\mathrm{Gal}}
\def\End{\mathrm{End}}
\def\Mor{\mathrm{Mor}}
\def\Aut{\mathrm{Aut}}
 \def\PB{\mathrm{PB}}
 \def\cA{\mathcal A}
\def\cB{\mathcal B}
\def\cC{\mathcal C}
\def\cD{\mathcal D}
\def\cE{\mathcal E}
\def\cF{\mathcal F}
\def\cG{\mathcal G}
\def\cH{\mathcal H}
\def\cJ{\mathcal J}
\def\cI{\mathcal I}
\def\cK{\mathcal K}
 \def\cL{\mathcal L}
\def\cM{\mathcal M}
\def\cN{\mathcal N}
 \def\cO{\mathcal O}
\def\cP{\mathcal P}
\def\cQ{\mathcal Q}
\def\cR{\mathcal R}
\def\cS{\mathcal S}
\def\cT{\mathcal T}
\def\cU{\mathcal U}
\def\cV{\mathcal V}
 \def\cW{\mathcal W}
\def\cX{\mathcal X}
 \def\cY{\mathcal Y}
 \def\cZ{\mathcal Z}
\def\0{{\ov 0}}
 \def\1{{\ov 1}}
 
 \def\frA{\mathfrak A}
 \def\frB{\mathfrak B}
\def\frC{\mathfrak C}
\def\frD{\mathfrak D}
\def\frE{\mathfrak E}
\def\frF{\mathfrak F}
\def\frG{\mathfrak G}
\def\frH{\mathfrak H}
\def\frI{\mathfrak I}
 \def\frJ{\mathfrak J}
 \def\frK{\mathfrak K}
 \def\frL{\mathfrak L}
\def\frM{\mathfrak M}
 \def\frN{\mathfrak N} \def\frO{\mathfrak O} \def\frP{\mathfrak P} \def\frQ{\mathfrak Q} \def\frR{\mathfrak R}
 \def\frS{\mathfrak S} \def\frT{\mathfrak T} \def\frU{\mathfrak U} \def\frV{\mathfrak V} \def\frW{\mathfrak W}
 \def\frX{\mathfrak X} \def\frY{\mathfrak Y} \def\frZ{\mathfrak Z} \def\fra{\mathfrak a} \def\frb{\mathfrak b}
 \def\frc{\mathfrak c} \def\frd{\mathfrak d} \def\fre{\mathfrak e} \def\frf{\mathfrak f} \def\frg{\mathfrak g}
 \def\frh{\mathfrak h} \def\fri{\mathfrak i} \def\frj{\mathfrak j} \def\frk{\mathfrak k} \def\frl{\mathfrak l}
 \def\frm{\mathfrak m} \def\frn{\mathfrak n} \def\fro{\mathfrak o} \def\frp{\mathfrak p} \def\frq{\mathfrak q}
 \def\frr{\mathfrak r} \def\frs{\mathfrak s} \def\frt{\mathfrak t} \def\fru{\mathfrak u} \def\frv{\mathfrak v}
 \def\frw{\mathfrak w} \def\frx{\mathfrak x} \def\fry{\mathfrak y} \def\frz{\mathfrak z} \def\frsp{\mathfrak{sp}}
 \def\bfa{\mathbf a} \def\bfb{\mathbf b} \def\bfc{\mathbf c} \def\bfd{\mathbf d} \def\bfe{\mathbf e} \def\bff{\mathbf f}
 \def\bfg{\mathbf g} \def\bfh{\mathbf h} \def\bfi{\mathbf i} \def\bfj{\mathbf j} \def\bfk{\mathbf k} \def\bfl{\mathbf l}
 \def\bfm{\mathbf m} \def\bfn{\mathbf n} \def\bfo{\mathbf o} \def\bfp{\mathbf p} \def\bfq{\mathbf q} \def\bfr{\mathbf r}
 \def\bfs{\mathbf s} \def\bft{\mathbf t} \def\bfu{\mathbf u} \def\bfv{\mathbf v} \def\bfw{\mathbf w} \def\bfx{\mathbf x}
 \def\bfy{\mathbf y} \def\bfz{\mathbf z} \def\bfA{\mathbf A} \def\bfB{\mathbf B} \def\bfC{\mathbf C} \def\bfD{\mathbf D}
 \def\bfE{\mathbf E} \def\bfF{\mathbf F} \def\bfG{\mathbf G} \def\bfH{\mathbf H} \def\bfI{\mathbf I} \def\bfJ{\mathbf J}
 \def\bfK{\mathbf K} \def\bfL{\mathbf L} \def\bfM{\mathbf M} \def\bfN{\mathbf N} \def\bfO{\mathbf O} \def\bfP{\mathbf P}
 \def\bfQ{\mathbf Q} \def\bfR{\mathbf R} \def\bfS{\mathbf S} \def\bfT{\mathbf T} \def\bfU{\mathbf U} \def\bfV{\mathbf V}
 \def\bfW{\mathbf W} \def\bfX{\mathbf X} \def\bfY{\mathbf Y} \def\bfZ{\mathbf Z} \def\bfw{\mathbf w}

 \def\R {{\mathbb R }} \def\C {{\mathbb C }} \def\Z{{\mathbb Z}} \def\H{{\mathbb H}} \def\K{{\mathbb K}}
 \def\N{{\mathbb N}} \def\Q{{\mathbb Q}} \def\A{{\mathbb A}} \def\T{\mathbb T} \def\P{\mathbb P} \def\G{\mathbb G}
 \def\bbA{\mathbb A} \def\bbB{\mathbb B} \def\bbD{\mathbb D} \def\bbE{\mathbb E} \def\bbF{\mathbb F} \def\bbG{\mathbb G}
 \def\bbI{\mathbb I} \def\bbJ{\mathbb J} \def\bbL{\mathbb L} \def\bbM{\mathbb M} \def\bbN{\mathbb N} \def\bbO{\mathbb O}
 \def\bbP{\mathbb P} \def\bbQ{\mathbb Q} \def\bbS{\mathbb S} \def\bbT{\mathbb T} \def\bbU{\mathbb U} \def\bbV{\mathbb V}
 \def\bbW{\mathbb W} \def\bbX{\mathbb X} \def\bbY{\mathbb Y} \def\kappa{\varkappa} \def\epsilon{\varepsilon}
 \def\phi{\varphi} \def\le{\leqslant} \def\ge{\geqslant}

\def\UU{\bbU}
\def\Mat{\mathrm{Mat}}
\def\tto{\rightrightarrows}

\def\Gr{\mathrm{Gr}}

\def\ch{\cosh}
\def\sh{\sinh}

\def\B{\bfB} 

\def\graph{\mathrm{graph}}

\def\gl{\mathfrak{gl}}

\def\la{\langle}
\def\ra{\rangle}

\def\V{\Updelta}
\def\ctg{\cot}

\def\e{\fre}

\begin{center}
	\bf\Large
 Projectors separating  spectra
 \\
  for $L^2$ on pseudounitary groups $\U(p,q)$
  
  \medspace
  
  \sc Yury A. Neretin%
  \footnote{Supported by the grants FWF, P25142, P28421}
  
\end{center}

{\small The spectrum of $L^2$ on a pseudo-unitary group $\U(p,q)$ (we assume $p\ge q$)
	 naturally splits into  $q+1$ types.
	 We write explicitly orthogonal projectors in $L^2$
	to subspaces with uniform spectra (this is an old question formulated by Gelfand and Gindikin).
	We also write two finer separations of
	$L^2$. In the first case pieces are enumerated by $r=0$, 1, \dots, $q$
	and  representations of discrete series of $\U(p-r,q-r)$, where
	$r=0$, \dots, $q$. In the second case pieces are enumerated by all discrete parameters of 
	the tempered spectrum of $\U(p,q)$.}

\section{Formulas for the projectors}

\COUNTERS

{\bf \punct Problem of separation of spectra.} Recall a problem 
formulated in the paper  \cite{GG}  by I.~M.~Gelfand and S.~G.~Gindikin in  1977.
Consider a real semisimple Lie group $G$, the left-right action of $G\times G$
on $G$ and the corresponding regular  representation 
in $L^2(G)$ (the group is equipped with the Haar measure).
The spectrum of the regular representation splits in a natural way
into several pieces  (according the number of non-conjugate Cartan subgroups).
Therefore there is a natural decomposition
of $L^2(G)$ into a direct sum of subrepresentations with uniform spectra,
$$
L^2(G)=L_1\oplus \dots \oplus L_m.
$$
Respectively, we have a natural decomposition of the identity operator
$$
E=\Pi_1+\dots+\Pi_m,
$$
where $\Pi_j$ are orthogonal projectors to the subspaces $L_j$. 
There arises a question about  explicit descriptions of such decompositions.

In \cite{GG} there was considered the case $G=\SL(2,\R)$. The space $L^2\bigl(\SL(2,\R)\bigr)$
is a sum of highest weight representations, a sum of lowest weight representations, and a
direct integral over the continuous series. It appears that  the summands corresponding
to highest weight and lowest weight representations can be regarded as
certain Hardy spaces
$H^2$.

The same question about separation of spectra arises for $L^2$
on semi-simple pseudo-Riemannian symmetric spaces and for some other problems of non-commutative harmonic analysis 
(a natural splitting of spectrum to different pieces is a  usual phenomenon).

\sm

{\bf \punct Known results.}
a) Transparent descriptions of decompositions are known  for several problems related
to $\SL(2)$:

\sm

--- For $L^2(\SL(2,\R))$, see \cite{GG}, \cite{Gin2}, \cite{Gin3}, \cite{AU}.

\sm

--- Consider the homogeneous space $\SL(2,\R)/H$, where $H$ is the diagonal
subgroup, this space can be identified with one-sheeted hyperboloid. The separation 
of spectrum in $L^2$  was discussed 
in \cite{Mol1}, \cite{Ner-rest}.

\sm

--- The space $L^2\bigl(\SL(2,\C)/\SL(2,\R)\bigr)$ was considered in \cite{FL}.

\sm

b) G.~I.~Olshanski \cite{Ols1} proposed a way of splitting of holomorphic series using 
non-commutative 'Hardy spaces', this approach was used in several works, see, e.g.,
\cite{HO}, \cite{BS}, \cite{Kof}, \cite{KM}.

Also, boundary values of holomorphic functions 
allow  to split off a part of a mostly continuous series for
$L^2$ on some pseudo-Riemannian symmetric spaces, \cite{GK}.

\sm

c) S.~G.~Gindikin \cite{Gin1} and V.~F.~Molchanov \cite{Mol2} in different ways solved a problem for multi-dimensional
hyperboloids $\OO(p,q)/\OO(p,q-1)$; this covers also all $\SL(2)$-cases mentioned above.

\sm

d) In \cite{Ner-compl} there was proposed a way to split summands of complementary series
using trace theorems,
see more in \cite{NO}, \cite{Ner-separ}.

\sm

However up  to now
explicit separations of spectra on groups  remain to be unknown except $\SL(2,\R)$.

In the present paper we  obtain such description for $L^2$ on pseudounitary
groups $\U(p,q)$.

\sm

{\bf \punct Pseudounitary groups and principal series.\label{ss:principal-series}}
Let $p\ge q>0$. We consider the space
$\C^p\oplus \C^q$ equipped with an Hermitian form with matrix 
$\begin{pmatrix}
1_p&\\ &-1_q
\end{pmatrix}$  (absent elements of matrices are 0 on default, $1_j$ denotes a unit matrix of
size $j$). 
The group 
$$
G:=
\U(p,q)
$$ consists of matrices $g$ preserving this form, i.e.,
\begin{equation}
g \begin{pmatrix}
1_p&\\ &-1_q
\end{pmatrix} g^*=\begin{pmatrix}
1_p&\\ &-1_q
\end{pmatrix}.
\label{eq:def-U}
\end{equation}
Consider  the {\it left-right regular representation} $R$
of the group $\U(p,q)\times\U(p,q)$ in $L^2\bigl(\U(p,q)\bigr)$,
$$
\cR(h_1,h_2) f(g):=f(h_1^{-1}g h_2),
\qquad \text{$(h_1,h_2)\in  \U(p,q)\times\U(p,q)$, $g\in \U(p,q)$}
.
$$
Recall a decomposition of $\cR$ into an integral of irreducible representations.

Denote by $J_r$ the $r\times r$-matrix with units on the secondary diagonal
(other matrix elements are 0-s).
For a given $r=0$, 1, \dots, $q$,
consider an Hermitian form  determined by a matrix
$$
I_r:=
\begin{pmatrix}
0& & & J_r\\
&1_{p-r}& &\\
&& 1_{q-r}&\\
J_r & & & 0
\end{pmatrix}
$$
For different $r$ these forms are equivalent.
Therefore 
the group of all  matrices $g$ satisfying 
$$
gI_r g^*=
I_r
.
$$
 is isomorphic to $\U(p,q)$.
 In this model, we consider subgroup $P_r\subset \U(p,q)$ of all block upper-triangular matrices $h\in \U(p,q)$  of size 
$$
\underbrace{1+\dots+1}_{\text{$r$ times}}+(p+q-2r)+\underbrace{1+\dots+1}_{\text{$r$ times}}
$$
having the form
$$
h=
\begin{pmatrix}
\ov \zeta_1^{\,-1}&*&\dots&*&*&* &\dots & * & *\\
&\ov \zeta_2^{\,-1}&\dots&*&*&* &\dots & * & *\\
&&\ddots&*&*&* &\dots & * & *\\
&&& \ov \zeta_r^{\,-1}&*&* &\dots & * & *\\
&&& & Z&* &\dots & * & *\\
&&& & & \zeta_r & \dots & * & *\\
&&& & &  & \ddots & * & *\\
&&& & &  &  & \zeta_2 & *\\
&&& & &  &  &  & \zeta_1
\end{pmatrix}.
$$
Here $Z\in \U(p-r,q-r)$. We consider a representation $\mu$ of $P_r$
given by
\begin{equation}
\mu_{A;c,m,\rho}(h)=\prod_{j=1}^r |\zeta_j|^{i\rho_j} \zeta_j^{m_j}
\cdot \tau_{A;c}(Z),
\label{eq:rep-parabolic}
\end{equation}
where $m_j\in \Z$, $\rho_j\in \R$, and $\tau_{A;c}$ is an irreducible representation
of $\U(p-r,q-r)$ of  discrete series%
\footnote{By the definition, a representation of a reductive 
group is contained in a {\it discrete series} if it is contained in $L^2$ on the group.}.
Below in Subsect. \ref{ss:characters} we will explain
the meaning of the parameters $(A;c)$.
Until this, we can understand $\tau_{A;c}$ as a symbol denoting an arbitrary representation of a 
discrete series of
$\U(p-r,q-r)$.
We consider  representations
$T_{A;c,m,\rho}$  of $\U(p,q)$
unitary induced (see, e.g, \cite{BR}, \S16.1) from  representations $\mu_{A;c,m,\rho}$. For $\rho$ being in a
 general position they are irreducible (see the Harish-Chandra completeness theorem,
 \cite{Kna}, Theorem 14.31). Thus we get $q+1$ family $\frA_r$ of representations numerated
 by $r=0$, 1, \dots, $q$. The regular left-right representation of $\U(p,q)\times \U(p,q)$
 admits a decomposition in a multiplicity free direct integral of the form
 \begin{equation}
 L^2\bigl(\U(p,q)\bigr)\simeq \bigoplus_{r=0}^q \int_{T_{A;c,m,\rho}\in \frA_r}
 (T_{A;c,m,\rho})^*\otimes T_{A;c,m,\rho}\,d\cP_r(A;c,m,\rho),
 \label{eq:decomposition-u}
 \end{equation}
 where $d\cP_r(A;c,m,\rho)$ is the Plancherel measure.

 \sm
 
 {\bf \punct Purpose of the paper.}
Our main purpose is to write projectors $\Pi_r$ corresponding to the orthogonal
 decomposition $\oplus_{r=0}^q$ in (\ref{eq:decomposition-u}).
 The formula is given in Theorem 1 at the end of this section.

 We also consider two  finer decompositions. 
 First, we fix a representation $\tau_{A;c}$ of a discrete series
 of $\U(p-r,q-r)$ 
 and consider in (\ref{eq:decomposition-u}) the integral of all representations
 having fixed $(A;c)$. Secondly, we fix $(A;c,m)$.
  Formulas for the corresponding orthogonal projectors are given in Theorems 2-3 in Section 2.
  Notice that for $r=0$ these formulas must coincide with characters of discrete series,
  general formulas have a similar degree of complexity (not too simple).

The problems are reduced to an integration of characters as functions
of parameters  with respect to the Plancherel measure%
\footnote{The Plancherel measure is supperted by a space with discrete and continuous coordinates, below we 
prefer to say {\it'summation' of characters}.}.
Characters of representations of real semisimple groups and the Plancherel formula
were obtained
by Harish-Chandra \cite{Harish-discrete}, \cite{Harish-Plancherel}. His formulas 
contain some undetermined constants, for more explicit formulas, see \cite{HW}.

We obtain formulas for the projectors $\Pi_r$
as a simple byproduct of Takeshi Hirai's \cite{Hir} derivation of the Plancherel
formula for $\U(p,q)$.  
 
 \sm
 
 {\bf \punct Cartan subgroups.} We realize $\U(p,q)$   as (\ref{eq:def-U}).
 Cartan subgroups 
 $H_k$, where $k=0$, \dots, $q$,
 are defined in the following way. First,  define a subgroup $H_k^+$ consisting of
 matrices
 {\small
 $$
 \begin{pmatrix}
 1_{p-k}&&& &&& &&& 0\\
 &\ch t_k & &&& &&&\sh t_k&\\
 &&\ch t_{k-1} && &&&\sh t_{k-1}&&\\
 && &\ddots&&&{}_{\text{\bf.}}  {\cdot}^{\text{\bf.}}&&&\\
 &&&&\ch t_1&\sh t_1&&&&\\
 &&&&\sh t_1&\ch t_1&&&&\\
  && & {}_{\text{\bf.}}  {\cdot}^{\text{\bf.}}  && & \ddots  && &\\
   &&\sh t_{k-1} && &&&\ch t_{k-1}&&\\   
   &\sh t_k & &&& &&&\ch t_k&\\
   0&&& &&& &&& 1_{q-k}
 \end{pmatrix}
 $$}
 Next, we define 
a subgroup $H_k^-$ consisting of diagonal matrices with entries
\begin{multline*}
\text{$e^{i\phi_1}$, $e^{i\phi_2}$, \dots, $e^{i\phi_{p-k}}$, $e^{i\theta_k}$,
$e^{i\theta_{k-1}}$, \dots, $e^{i\theta_1}$,}
\\
\text{	$e^{i\theta_1}$,
\dots, $e^{i\theta_{k-1}}$, $e^{i\theta_k}$, $e^{i\psi_{q-k}}$,
\dots, $e^{i\psi_2}$, $e^{i\psi_1}$.
	}
\end{multline*}
Here $t_\gamma\in \R$, $\phi_\alpha$, $\psi_\beta$, $\theta_\gamma\in\R/2\pi\Z$.
We set
$$
H_k:=H_k^+\cdot H_k^-.
$$
Denote
$$
z_j:=t+i\theta_j.
$$
The eigenvalues of $h\in H_k$ are
\begin{multline}
\text{$e^{i\phi_1}$, $e^{i\phi_2}$, \dots, $e^{i\phi_{p-k}}$, $e^{z_1}$,  $e^{z_2}$,\dots, $e^{z_k}$,}
\\
\text{$e^{i\psi_1}$, $e^{i\psi_2}$, \dots, $e^{i\psi_{q-k}}$, $e^{-\ov z_1}$,  $e^{-\ov z_2}$,\dots, $e^{-\ov z_k}$.}
\label{eq:eigenvalues}
\end{multline}
It is convenient to use two notations for systems of coordinates on $H_k$, the first is
$
\phi_\alpha, \psi_\beta, \theta_\gamma, t_\gamma,
$
the second is
$
\phi_\alpha, \psi_\beta, z_\gamma
$.

Define the canonical Lebesgue measure on $H_k$ by
$$
d_kh=\prod_{\alpha=1}^{p-k} d\phi_\alpha \prod_{\beta=1}^{q-k} d\psi_\beta
\prod_{\gamma=1}^k dt_\gamma\,d\theta_\gamma.
$$

The  {\it Weyl group} $W_k$ corresponding to a Cartan subgroup $H_k$ is 
$$
W_k\simeq S_{p-k}\times S_{q-k} \times (S_k\ltimes \Z_2^k),
$$
 the symmetric group $S_{p-k}$ acts on $H_k$ by permutations of coordinates
 $\phi_\alpha$, the group $S_{q-k}$ by permutations of $\psi_\beta$,
 the $S_k$ by permutations of pairs $(t_\gamma,\theta_\gamma)$, and $\Z_2^k$
 is generated
 by $k$ reflections 
 $$
 R_\gamma: (t_1,\dots,t_{\gamma-1}, t_\gamma, t_{\gamma+1}, \dots, t_k)\mapsto
 (t_1,\dots,t_{\gamma-1}, -t_\gamma, t_{\gamma+1}, \dots, t_k)
 $$
 (over coordinates $\phi_\alpha$, $\psi_\beta$, $\theta_\gamma$ on $H_k$ remain fixed).

 We say that a function
 $f$ on $H_k$ is {\it $\epsilon_k$-symmetric} if it is invariant with respect
 to the subgroups $S_{p-k}$, $S_{q-k}$, $S_k$ and changes sign
 under each reflections $R_\gamma$. We say that a function $f$ is
 {\it $\epsilon_k$-skew-symmetric} if it is skew-symmetric with respect to
$S_{p-k}\times S_{q-k}$ and invariant with respect to $S_k\ltimes \Z_2$. 
 
\sm

{\bf\punct The Vandermonde expression.} We denote by $\V(y)$ the Vandermonde expression
$$
\V(y)=\V(y_1,\dots,y_n)= \prod_{1\le j <l\le n}^n (y_j-y_l).
$$

Denote the eigenvalues (\ref{eq:eigenvalues})
 of a matrix $h\in H_k$ by $e^{x_1}$, \dots, $e^{x_n}$, and set
$$
\V_k(h):=\V(e^{x_1},\dots,e^{x_n}).
$$ 

Next, consider the following differential operators
$X_1$, \dots, $X_n$ on $H_k$:
\begin{multline}
\frac\partial{i\partial \phi_1},\dots, \frac\partial{i\partial \phi_{p-k}},\,
\frac12\Bigl(\frac \partial {\partial t_1}+ \frac{\partial}{i \partial\theta_1}\Bigr),\dots,
\frac12\Bigl(\frac \partial {\partial t_k}+ \frac{\partial}{i \partial\theta_k}\Bigr)
\\
\frac\partial{i\partial \psi_1},\dots, \frac\partial{i\partial \psi_{q-k}},\,
\frac12\Bigl(-\frac \partial {\partial t_1}+ \frac{\partial}{i \partial\theta_1}\Bigr),\dots,
\frac12\Bigl(-\frac \partial {\partial t_k}+ \frac{\partial}{i \partial\theta_k}\Bigr).
\label{eq:XXX}
\end{multline}
We set%
\footnote{Such operators were introduced 
	Gelfand and Naimark \cite{GN} for complex classical groups and later were  used by
	Harish-Chandra in a more general context.}
$$
\V_k(\partial)=\V(X_1,\dots,X_n).
$$

\sm

{\bf\punct Average operator.%
\label{ss:average}} 
Denote by $C_c^\infty(G)$ the space of compactly supported smooth functions on $\U(p,q)$.
Recall that a space $G/H_k$ admits a unique up to a scalar factor
$G$-invariant measure (since both $G$, $H$ are unimodular, see, e.g., \cite{BR},\S 4.3).
For any $f\in C_c^\infty(G)$ we assign a function ({\it Harish-Chandra transform} of
$f$) on $H_k$ by
$$
I_k f(h)=\int_{G/H_k} f(yhy^{-1}) dy, \qquad h\in H_k,\quad y\in G.
$$
Notice that for $y\in G$ the expression $yhy^{-1}$ depends only on a coset
of $G$ with respect to $H_k$, therefore actually we  have an integration over $G/H_k$.
By the definition, a function $I_k f$ is invariant with respect 
to the group $W_k$.

Under a certain normalization of
the Haar measure $dg$ on $G$ and invariant measures on $H_k$ we have a {\it Weyl integration formula}
$$
\int_G f(g)\,dg =  \sum_{k=0}^q \omega_k \int_{H_k} I_k f(h)\, |\V(h)|^2 \, d_k(h),
$$
where 
$$
\omega_k :=\frac 1{\# W_k}=\frac 1{(p-k)!\,(q-k)!\, k!\,2^k}
$$
(this  follows from the usual arguments establishing the Weyl integration formula).

Denote by $H_k'\subset H_k$ the subset $\V_k(h)\ne 0$.
By $H_k^\circ\subset H_k$ we denote the complement to the 
union of all hypersurfaces $\phi_j=\psi_l$.
For any $f\in C_c^\infty(G)$ we define a function 
on $H_k'$
 by
\begin{equation}
\Xi_k f(h):= \prod \sgn (t_k) \cdot \ov{\V(h)} I_k f(h).
\label{eq:Xi-k}
\end{equation}

A function $\Phi=\Xi_k f$ satisfies the following properties
(see \cite{War}, Corollary 8.5.1.2).

\sm

A) $\Phi$ is compactly supported;

\sm

B)  $\Phi$ is $\epsilon_k$-skew-symmetric;

\sm

C) $\Phi$  is $C^\infty$-smooth on each component $U$ of 
$H_k^\circ$ 
and all partial derivatives are bounded.
 
\sm

Moreover, an a operator $\Xi_k$ is bounded 
in a natural sense, i.e. a convergence of a sequence
$f_j\in C_c^\infty(G)$ 
implies uniform convergence of all partial derivatives 
of $\Xi_k f$ on $H_k^\circ$.

\sm

A collection of functions $(\Xi_0 f, \dots, \Xi_q f)$
satisfies some gluing conditions on hypersurfaces $\phi_\alpha=\psi_\beta$
and $t_\gamma=0$, see Lemma 2.2 of \cite{Hir}.
These conditions are an important element of the story about characters and 
the Plancherel formula, but below
we do not use them explicitly (these conditions are hidden in an integration by parts 
in formula (\ref{eq:by-parts})).

Next, consider a function $F:=\V_k(\partial) \Xi_k f$.
It satisfies the following conditions (Harish-Chandra \cite{Harish}, Lemma 40, in \cite{Hir}
it is formulated in a beginning of Sect.3, see also \cite{War}, Corollary 8.5.1.5):

\sm

A$^\circ$) $F$ is compactly supported;

\sm

B$^\circ$)  $F$ is $\epsilon_k$-symmetric;

\sm

C$^\circ$) $F$ admits a $C^\infty$-smooth extension to the closure 
of each component  of $H_k^\circ$ ;

\sm

D$^\circ$) $F$ admits a continuous extension to the whole $H_k$.

\sm


\sm

{\bf\punct Formula for projectors.}
Consider a distribution $\chi$ on $G$ invariant with respect to conjugations.
It determines a convolution operator
$Q_\chi:C^\infty_c(G)\to C^\infty(G)$ by%
\footnote{By $\lla f,\chi\rra_{(L)}$ we denote a pairing of a test function and
	a distribution on a manifold $L$.}
$$
Q_\chi f(h)=\lla  f(gh), \chi(g) \rra_{(G)}
,$$
where brackets denote a pairing of a test function and a distribution.

\begin{ttheorem}
 For any $r=0$, \dots, $q$
 the projector $\Pi_r$ is a convolution operator determined by
 the following distribution $\Theta_r$:
 \begin{multline}
 M_* \lla f,\Theta_r\rra_{(G)}=
(-1)^{n(n-1)/2+pq+qr+r(r-1)/2}  \sum_{k=r}^q 
 \frac{ 2^{n-2k}\pi^{n-k} k! }{(k-r)!\, r!}
 \times \\ \times
 \omega_k\Bigl\la\!\!\Bigl\la \V_k(\partial)\, \Xi_k \,f(\phi,\psi,\theta,t),
\\\prod_{\alpha=1}^{p-k} \delta(\phi_\alpha)
\prod_{\beta=1}^{q-k} \delta(\psi_\beta) 
 \prod_{\gamma=1}^k \bigl(\coth (t_\gamma/2)\,\delta(\theta_\gamma)+
 \tanh(t_\gamma/2)\delta(\theta_\gamma-\pi)\bigr)\Bigr\ra\!\!\Bigr\ra_{(H_k)},
 \label{eq:th1}
  \end{multline}
  where $\delta(\cdot)$ denotes the delta-function and $M_*$ is a constant.
\end{ttheorem}

{\sc Remark.} A formula
\begin{equation}
\lla F,\coth(t/2)\rra_{(\R)}:=
\mathrm{p.v}
\int_{\R}\coth(t/2) F(t)\,dt
\label{eq:coth-pv}
\end{equation}
determines a distribution on $C_c^\infty(\R)$. However test
functions in (\ref{eq:th1}) are odd with respect to the variables $T_\gamma$, for odd smooth $F(t)$ the integrand
in (\ref{eq:coth-pv}) is smooth at 0.

\sm 

{\sc Remark.}
In particular, the projector to the most continuous series (i.e., 
$r=q$) is determined by the distribution
\begin{multline*}
 M_* \lla f,\Theta_q\rra_{(G)}=
 (-1)^{p(p-1)/2+q^2}   
 2^{p-q}\pi^{p} \omega_q
 \times \\ \times
\Bigl\la\!\!\Bigl\la \V_q(\partial)\, \Xi_q \,f(\phi,\psi,\theta,t),
\\
 \prod_{\alpha=1}^{p-q} \delta(\phi_\alpha) 
 \prod_{\gamma=1}^q \bigl(\coth (t_\gamma/2)\,\delta(\theta_\gamma)+
 \tanh(t_\gamma/2)\delta(\theta_\gamma-\pi)\bigr)\Bigr\ra\!\!\Bigl\ra_{(H_q)},
 \end{multline*}

{\bf\punct Further structure of the paper.}
Section 2 contains preliminaries from Hirai \cite{Hir}.
Theorem 1 is proved in Section 3. 
In Section 4 we write formula for projectors
determining finer orthogonal decompositions of 
$L^2\bigl(\U(p,q)\bigr)$.

\section{The Plancherel formula. Preliminaries}

\COUNTERS

\def\zcircle{\epsfbox{char.4}}
\def\zotimes{\epsfbox{char.5}}
\def\zblackcircle{\epsfbox{char.6}}
\def\zsquare{\epsfbox{char.7}}
\def\zboxtimes{\epsfbox{char.8}}
\def\zblacksquare{\epsfbox{char.9}}

\sm

Here we present the formula for characters and the Plancherel formula
from \cite{Hir}.

\sm

{\bf \punct Formula for characters.%
\label{ss:characters}}
Recall that a character $\pi$
of a unitary representation $T$ of a unimodular Lie group $G$ is a distribution on $G$
defined by
$$
\lla f,\pi\rra_{G)}=
\tr T(f):=\tr \int_G f(g)\,T(g)\,dg
\qquad \text{for all $f\in C_c^\infty(G)$.}
$$
According Harish-Chandra, a character of an irreducible representation of a reductive Lie group
is a locally integrable function.
Here we present a formula for characters of representations of $\U(p,q)$ 
from \cite{Hir}, Section 1.

Fix $r=0$, $1$, \dots, $q$.
Consider three collections of parameters
\begin{align}
&\text{$c_1>c_2>\dots>c_{n-2r}$, \quad where $c_\alpha\in \Z$};
\label{eq:c}
\\
&\text{$m_1$, \dots, $m_r$,\quad where $m_j\in \Z$};
\\
&\text{$\rho_1>\rho_2>\dots>\rho_r>0$, where $\rho_j\in\R$.}
\label{eq:rho}
\end{align}
We  denote
$$
d_j=\tfrac 12 (m_\gamma+i\rho_\gamma), 
$$
and use an alternative notation (a {\it signature}) for the same collection of parameters
\begin{equation}
(c,d):=
(c_1, c_2 \dots, c_{n-2r}, d_1,\ov d_1, d_2,\ov d_2,\dots, d_r, \ov d_r). 
\label{eq:cd}
\end{equation}

Next, we split the set $\{1,2,\dots,n-2r\}$ into two disjoint subsets
$A:=\{a_1,\dots,a_{p-r}\}$ and $B:=\{b_1,\dots,b_{q-r}\}$.
Data $(A;c,d)$ determine  a character
$\kappa_{A;c,d}$ of the group $\U(p,q)$, it is defined in this subsection.

Let $j$, $l\in \Z$, $z=t+i\theta$ .
We set
\begin{equation}
 \xi_c(z;j,l):=\sgn(j-l) \exp\bigl\{-|c_j-c_l|\, |t|+(c_j+c_l) i\theta \bigr\}.
\end{equation}

\begin{figure}
 $$
\epsfbox{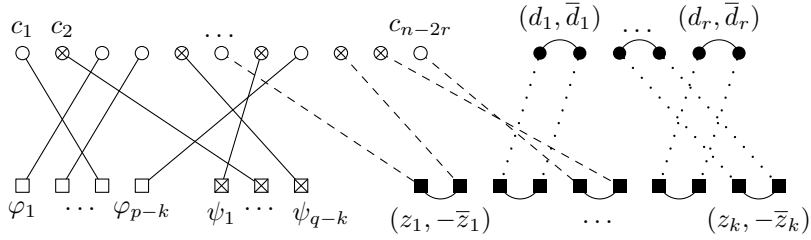}
$$
\caption{Diagrams, which enumerate summands in the formula for characters.\label{fig:diagram}}
\end{figure}

For  given $r$ and $k\ge r$ we 
consider all possible diagrams $\sigma$ of the form given on Fig.\ref{fig:diagram}.

The elements of the upper row, $\zcircle$-s, $\zotimes$-s, $\zblackcircle$-s
correspond to elements of a  signature (\ref{eq:cd}).
More precisely, $\zcircle$-s correspond to $c_{a_j}$, where
$a_j\in A$, $\zotimes$-s correspond to $c_{b_l}$, where
$b_l\in B$, and $\zblackcircle$-s correspond to $d_1$, $\ov d_1$, \dots,
$d_r$, $\ov d_r$.
For a clarity, we connect $d_\gamma$ and $\ov d_\gamma$ by an arc.

The elements of the lower row, $\zsquare$-s, $\zboxtimes$-s, and $\zblacksquare$-s 
correspond to the coordinates
\begin{equation}
\phi_1,\dots,\phi_{p-k}, \psi_1,\dots,\psi_{q-k}, z_1,-\ov z_1, z_2, -\ov z_2, 
\dots, z_k,-\ov z_k.
\label{eq:coordinates}
\end{equation}
Namely, $\zsquare$-s correspond to $\phi_\alpha$, $\zboxtimes$-s to $\psi_\beta$,
and $\zblacksquare$-s to $z_\gamma$, $-\ov z_\gamma$. We connect 
each pair $z_\gamma$ and $-\ov z_\gamma$ by an arc.

We connect elements of the upper row with elements of the lower row by arcs
(each element is an end of a unique arc). Each diagram $\sigma$
establishes a one-to-one correspondence
between elements of rows (\ref{eq:cd}) and of rows (\ref{eq:coordinates}).

We allow only diagrams that are unions of pieces of 4 types a)-d) presented on 
Fig. \ref{fig:pieces}:
\begin{figure}
$$
\epsfbox{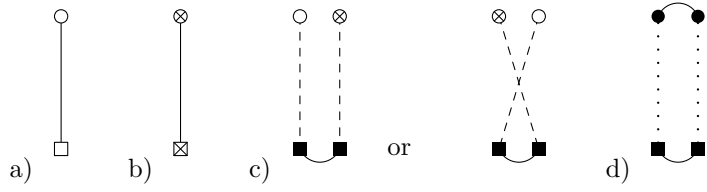}
$$ 
\caption{Possible pieces of a diagram.\label{fig:pieces}}
\end{figure}

a) Arcs $\zcircle$---$\zsquare$ or  $c_{a_j}$---$\phi_\alpha$, where
$a_j\in A$.

\sm

b) Arcs $\zotimes$---$\zboxtimes$ or $c_{b_l}$---$\psi_\beta$, 
where $b_l\in B$.

\sm

c) Chains $\zcircle$---$\zblacksquare$---$\zblacksquare$---$\zotimes$ 
or $c_{a_j}$---$z_\gamma$---$(-\ov z_\gamma)$---$c_{b_l}$.
Un particular, this means  that a {\it left} $\zblacksquare$ is connected with  $\zcircle$
and a {\it right} $\zblacksquare$ is connected with $\zotimes$.

\sm

d) Cycles --$\zblackcircle$---$\zblacksquare$---$\zblacksquare$---$\zblackcircle$--
or --$d_s$---$z_\gamma$---$(-\ov z_\gamma)$---$\ov d_s$--.
Notice that a {\it left} $\zblacksquare$ is connected with {\it left}
$\zblackcircle$.

\sm

We use the following notation:

\sm

--- in the case a) we write $[c_{a_j}, \phi_\alpha]\in \sigma$;

\sm

--- in the case b): $[c_{b_l},\psi_\beta]\in\sigma$;

\sm

--- in the case c): $[c_{a_j}, z_\gamma, c_{b_l}]\in \sigma$;

\sm

--- in the case d): $[[d_s,z_\gamma]]\in\sigma$.

\sm

Denote by $\Omega(A|r,k)$ the set of all admissible diagrams $\sigma$.
Recall that $\sigma\in \Omega(A|r,k)$ determines a substitution
of a set $\{1,\dots,n\}$, in particular it has a well-defined
sign $(-1)^\sigma$. 

For a fixed signature (\ref{eq:cd}) and a fixed $A$ we define functions 
$\kappa^k=\kappa^k_{A;c,d}$ on $H_k$ in the following way:
$$
\kappa^k_{A;c,d}=0, \qquad \text{for $k<r$}.
$$
For $k\ge r$ we set
\begin{multline}
 \kappa^k_{A;c,d} (h)=
 (-1)^{k(k+1)/2+pq-r(k+q)}
 \sum_{\sigma\in \Omega(A|r,k)}
 (-1)^\sigma \times  \\\times 
 \prod_{[c_{a_j}, \phi_\alpha]\in \sigma} e^{i c_{a_j}\phi_\alpha}
 \cdot 
 \prod_{[c_{b_l},\psi_\beta]\in\sigma} e^{i c_{b_l} \psi_\beta}
 \times\\ \times
 \prod_{[c_{a_j}, z_\gamma, c_{b_l}]\in \sigma}
 \xi_c(z_\gamma;a_j,b_l)
 \cdot
 \prod_{[[d_s,z_\gamma]]\in\sigma} e^{im_s \theta_\gamma}
 \bigl( e^{i\rho_s t_\gamma}+e^{-i\rho_s t_\gamma}\bigr)
 .
 \label{eq:kappa-k}
\end{multline}
 The functions $\kappa^k_{A;c,d}$ are  $\epsilon_k$-skew-symmetric.

There exists a unitary representation $T_{A;c,d}$ of $G=\U(p,q)$
such that for any $f\in C_c^\infty(G)$, we have
$$
\tr  T_{A;c,d}(f)=\lla f(g), \pi_{A;c,d}\rra_{(G)}=
\sum_{k=0}^q\omega_k \int_{H_k} \Xi_k f(h)\,\kappa^k_{A;c,d}(h)\,dh.
$$

\sm

Moreover, $T_{A;c,d}$ are the representations defined in Subsect. \ref{ss:principal-series}.

If $r=0$ (i.e., the parameters $d$ are absent), then $T_{A;c}$ is a representation of 
$\U(p,q)$ of discrete series.

Let $r>0$. Then
$\kappa_{A;c,\varnothing}$ determines a representation of a discrete series $\tau_{A;c}$ of $\U(p-r,q-r)$.
The representation $T_{A;c,d}$ is induced from the representation (\ref{eq:rep-parabolic})
of the parabolic $P_r$.

\sm

{\bf \punct The Plancherel formula.} See Hirai \cite{Hir}, Theorem 3.
For $m\in \Z$, we set
$$
\e_m(\rho):=
\begin{cases}
 -(i/2) \coth (\pi\rho/2), \qquad \text{if $m$ is even};
 \\
 -(i/2) \tanh (\pi\rho/2), \qquad \text{if $m$ is odd}.
\end{cases}
$$
We also define the Vandermonde expression in the parameters $(c,d)$,
$$
\V_r(c,d)=\V(c_1,\dots,c_{p-r},d_1, \dots,d_r, c_{p-r+1},\dots, c_{n-2r},\ov d_1,\dots,\ov d_r)
.
$$

The Plancherel formula for $\U(p,q)$ is given by
\begin{multline}
M_* f(e)=\sum_{r=0}^q 
\biggl\{
\sum_{c_1>c_2>\dots>c_{n-2r}}\,\,
\sum_{m_1, \dots, m_r\in \Z}
\\
\int\limits_{\rho_1>\rho_2>\dots>\rho_r>0}
\Bigl(\sum_A \tr T_{A;c,d} (f)\Bigr)\, \V_r(c,d)\prod_s \e_{m_{s}}(\rho_s)\,d\rho_1\dots d\rho_r.
\biggr\}
\label{eq:plancherel}
,\end{multline}
where $M_*$ is a constant.

Denote by $F_1*F_2$ the convolution of functions on $G$.
For $F\in C_c^\infty(G)$ denote by $F^\star$ the function
$F^\star(g)=\ov{F(g^{-1})}$. For any unitary representation
$T$ of $G$ we have
$$
T(F*F^\star)=T(F) \,T(F)^*\ge 0.
$$ 
Therefore for $f=F*F^\circ$ the integrand in the right hand side of
(\ref{eq:plancherel}) is positive. By  polarization
arguments this implies absolute convergence
of the integral and the series in (\ref{eq:plancherel})
for functions of the form 
$f=F*F'$,  hence the absolute convergence holds
 on the subspace consisting of functions
 \begin{equation}
f= \sum_{j=1}^N F_j*F_j',\qquad\text{where $F_j$, $F_j'\in C_c^\infty(G)$}.
 \label{eq:subspace}
 \end{equation}
 This subspace is dense in $C^\infty_c(G)$ and invariant with respect to left and right shifts.
  For arbitrary $f\in C^\infty_c(G)$ the identity
(\ref{eq:plancherel}) holds if to understand the right-hand side in the sense of
a successive integration as below.

\sm


\sm

\section{Evaluation of the projectors $\Pi_r$}

\COUNTERS

Here we prove Theorem 1.

\sm

{\bf \punct Preliminary remarks.} 
Denote by $\wh G$ the set of all possible parameters
$\lambda=(A;c,d)$, see (\ref{eq:c})--(\ref{eq:rho}), so $\wh G$ consists
of pieces enumerated by $r=0$, 1, \dots, $q$, and each 
piece is a product of a discrete set and a 
simplicial cone $\rho_1>\dots>\rho_r>0$. We equip the discrete set 
with the counting measure and the simplicial cone with the Lebesgue
measure, so we get a sigma-finite measure on $\wh G$. 
 Denote  $d\cP(A;c,d)=d\cP(\lambda)$  the  measure
on $\wh G$
with the positive density
$\V_r(c,d)\prod_s \e_{m_{s}}(\rho_s)$.
Denote by $T_\lambda$ the irreducible representation
with parameter $\lambda$ and by $H_\lambda$ the space of the representation.

Next, consider the space $\cL^2(G)$ of functions $\Phi$ on $\wh G$,
which for each $\lambda$ assign a Hilbert--Schmidt operator
$\Phi(\lambda):H_\lambda\to H_\lambda$ and satisfy the condition
$$
\int_{\wh G} \tr
\bigl( \Phi(\lambda)\Phi(\lambda)^*\bigr)\,d\cP(\lambda)<\infty.
$$ 
This space is a Hilbert space with respect to
the inner product
$$
\la \Phi,\Psi\ra_{\cL^2(\wh G)}:= \int_{\wh G} \tr \Phi(\lambda)\Psi(\lambda)^*\,d\cP(\lambda).
$$
 For any 
$f\in L^2(G)\cap L^1(G)$ 
the formula
$$
T_\lambda(f)=\int_G f(g) T_\lambda(g)\,dg
$$
determines an element of $\cL^2(\wh G)$. Moreover
$$
\la f_1,f_2\ra_{L^2(G)}=\la T_\lambda(f_1),T_\lambda(f_2)\ra_{\cL^2(\wh G)}
$$
and the map $I:f\mapsto T_\lambda(f)$ (the {\it Fourier transform})
 extends to a unitary operator
$L^2(G)\to\cL^2(\wh G)$. The {\it inverse Fourier transform} is given by
$$
I^{-1} \Phi(g)=\int_{\wt G} \tr \bigl(\Phi(\lambda) T(g^{-1})\bigr)
\,d\cP(\lambda).
$$

Consider a subset $U\subset \wh G$ and  the subspace
$\cL^2(U)\subset \cL^2(\wh G)$ consisting of functions
supported by $U$. Let us write a formula for a projection
operator $\Pi_U$ in $L^2(G)$ corresponding to $\cL^2(U)$. According
the kernel theorem any bounded operator in $L^2(G)$ is determined by
a kernel, which is a distribution on $G\times G$. Let
$f$ has the form (\ref{eq:subspace}). Then
$$
\Pi_U f(g)=\int_{U} \tr\bigl( T_\lambda(f) T(g^{-1})\bigr)\,
d\cP(\lambda)
=\int_U \lla L_g f,\pi_\lambda\rra_{(G)}\,d\cP(\lambda),
$$ 
where $L_g f(h):=f(g^{-1}h)$
and the integral absolutely converges.

\sm

\sm

{\bf\punct Transformations of the Plancherel formula.%
\label{ss:transformations}}
 Denote by $W^{(r)}$ the group  of all transformations
of the set of all signatures (\ref{eq:cd}) generated by permutations 
of parameters $c$, permutations of parameters $d$
and reflections 
$
R^s
$
changing $d_s$ and $\ov d_s$. Clearly,
$$
W^{(r)}\simeq S_{n-2r}\times \bigl(S_r\ltimes \Z_2^r\bigr).
$$
We say that a function $F$
in variables $c$, $d$ is {\it $\epsilon^r$-symmetric} if it 
is invariant with respect to $S_{n-2r}$ and $S_r$
and changes a sign under each reflection $R^s$.
A function $F$ is {\it $\epsilon^r$-skew-symmetric}
if it is skew-symmetric with respect to $S_{n-2r}$
and invariant with respect to $S_r\ltimes \Z_2^r$.

We write the right-hand side of  (\ref{eq:plancherel}) as 
$$
\sum_{r=0}^q 
Z_r.
$$
Our purpose is to find summands $Z_r$.
Following \cite{Hir},
we define a sum
$$
\sum_A  \tr T_{A;c,d} (f)=
 \sum_k \omega_k \int_{H_k} \Xi_k f(h) \Bigl(\sum_A \kappa^k_{A;c,d}(h) \Bigr) \,dh.
$$
The expressions
$$
\wt \kappa_{c,d}^k:=
\sum_A \kappa^k_{A;c,d}
$$
have form (\ref{eq:kappa-k})
but the summation
$$
\sum_{\wt\Omega(r,k)} \bigl(\dots \bigr)
$$
\begin{figure}
$$
\epsfbox{char.3}
$$ 
\caption{Diagrams enumerating summands in $\wt\kappa^k_{c,d}$.%
\label{fig:diagram-2}}
\end{figure}
%
%
%
now is taken
over the set $\wt\Omega(A;r,k)$ of diagrams shown on Fig. \ref{fig:diagram-2},
namely we forget a difference between $\zcircle$-s and $\zotimes$-s
and allow to connect any $\phi_\alpha$, $\psi_\beta$ with an arbitrary $c_\tau$.
The number of elements
of $\wt\Omega(A|r,k)$ is
$$
\#\wt\Omega(A|r,k)=\frac{k!(n-2r)!}{(k-r)!}.
$$

We  defined $\kappa^k_{A;c,d}$ and $\wt\kappa^k_{c,d}$ under conditions
(\ref{eq:c}) and (\ref{eq:rho}). However, the expressions make sense for
arbitrary $c_1$, \dots, $c_{n-2r}\in \Z$ and $\rho_1$, \dots $\rho_r\in\R$.
Functions $\wt\kappa^k_{c,d}$ are $\epsilon^r$-skew-symmetric with respect to the parameters
 $c$, $d$, $\ov d$.

For any signature (\ref{eq:cd}) we define  functions $\eta^k_{c,d}$ on $H_k$ by
\begin{multline}
 \eta^k_{c,d}(h)=\sum_{\sigma\in \wt\Omega(r,k)}
  \prod_{[c_{a_j}, \phi_\alpha]\in \sigma} e^{i c_{a_j}\phi_\alpha}
 \cdot 
 \prod_{[c_{b_l},\psi_\beta]\in\sigma} e^{i c_{b_l} \psi_\beta}
 \times\\ \times
 \prod_{[c_{a_j}, z_\gamma, c_{b_l}]\in \sigma}
\Bigl(\sgn (t_\gamma) \exp\{-|c_{a_j}-c_{b_l}|\cdot |t_j|+ i(c_{a_j}+c_{b_l})\theta_j \}  \Bigr)
 \times\\ \times
 \prod_{[[d_s,z_\gamma]]\in\sigma}  e^{im_s\theta_\gamma} \bigl(e^{i\rho_s \theta_\gamma}- 
 e^{-i\rho_s \theta_\gamma}\bigr).
\label{eq:eta}
\end{multline}
This expression is $\epsilon_k$-symmetric as a function in variables
$\phi_\alpha$, $\psi_\beta$ and $\epsilon^r$- symmetric as a function
in parameters $c_a$, $d_s$, $\ov d_s$. It is easy to verify that
$$
\V_k(\partial) \eta^k_{c,d}(h)=(-1)^{pq+qr+r(r-1)/2} \V_k(c,d) \wt\kappa^k_{c,d}.
$$

{\sc Remark.}
We can also define  functions $\eta^k_{c,d|A}$ replacing a summation with respect
to $\wt\Omega(r,k)$ by a summation with respect to $\Omega(A|r,k)$.
These functions are $\epsilon_k$-symmetric in coordinates, but the $\epsilon^r$-symmetry with respect to
the parameters does not hold.
\hfill $\boxtimes$

\sm

Therefore, we can present $Z_r$ as
\begin{multline*}
 Z_r=(-1)^{pq+qr+r(r-1)/2}
 \sum_{c_1>c_2>\dots>c_{n-2r}}\,\,
\sum_{m_1, \dots, m_r\in \Z}
\\
\int\limits_{\rho_1>\rho_2>\dots>\rho_r>0}
\biggl( \sum_{k=r}^q \omega_k\int_{H_k}\Xi_k f(h)\cdot \V_k(\partial)\eta^k_{c,d}(h)
\,dh\biggr)  \prod_s \e_{m_{s}}(\rho_s)\,d\rho_1\dots d\rho_r.
\end{multline*}
Next, Theorem 2 of \cite{Hir} allows to integrate by parts:
\begin{multline}
 \sum_{k=r}^q \omega_k\int_{H_k}\Xi_k f(h)\cdot \V_k(\partial)\eta^k_{c,d}(h)\,dh
 =\\=
 (-1)^{n(n-1)/2}
 \sum_{k=r}^q \omega_k\int_{H_k}\V_k(\partial) \Xi_k f(h)\cdot \eta^k_{c,d}(h)
\,dh.
\label{eq:by-parts}
\end{multline}

{\sc Remark.}
This is a delicate point since functions $\Xi_k f(h)$ and $\eta^k_{c,d}(h)$
have singularities on hypersurfaces $\phi_\alpha=\psi_\beta$,
an integration by parts in each summand produces extra terms,
however in the sum $\sum_{k=r}^q$ all such terms cancel.
\hfill $\square$

\sm

{\bf \punct Transformations of distributions $Z_{r}$.}
Next, both $\V_k(\partial)\Xi_k f(h)$ and $\eta^k_{c,d}(h)$
are $\epsilon^r$-symmetric as functions of parameters $c$, $d$, $\ov d$,
Therefore, we can transform $Z_r$ to the form
\begin{multline*}
Z_r=
(-1)^{n(n-1)/2} (-1)^{pq+qr+r(r-1)/2}
\frac 1{\# W^{(r)}} 
 \sum_{c_1,\dots, c_{n-2r}\in \Z}\,\,
 \sum_{m_1, \dots, m_r\in \Z}
 \\
 \int\limits_{\rho_1, \dots,\rho_r\in \R}
 \biggl(\sum_{k=r}^q \omega_k \int_{H_k}\V_k(\partial)\Xi_k f(h)\cdot \eta^k_{c,d}(h)\,dh\biggr)
 \prod_s \e_{m_{s}}(\rho_s)\,d\rho_1\dots d\rho_r
\end{multline*}
(we also use the property $\e_m(\rho)=-\e_m(-\rho)$).

At least for $f$ of the form (\ref{eq:subspace}) 
this expression converges as an integral over $\wh G$ 
(the integrand is the expression in the big brackets). But we have also a finite summation
and an integration over $H_k$ and we have no reasons
to believe to the absolute convergence of the whole expression. For further manipulations we pass
to a successive integration and after this change the order of the successive
integration and a finite summation.
In this way, we transform $Z_r$
to the form
$$
Z_r=
\frac{(-1)^{n(n-1)/2} (-1)^{pq+qr+r(r-1)/2}}
 {\# W^{(r)}}
 \sum_{k=0}^q \omega_k
 Z_{r,k},
$$
where
\begin{multline}
Z_{r,k}= 
\sum_{c_1\in \Z} \dots \sum_{c_{n-2r}\in \Z}
\,\,
\sum_{m_1\in \Z} \int_{\rho_1\in\R}\sum_{m_2\in \Z} 
\int_{\rho_2\in\R}\dots
 \sum_{m_r\in \Z} \int_{\rho_r\in\R}
\\
\biggl(\int_{H_k}\V_k(\partial)\Xi_k f(h)\cdot \eta^k_{c,d}(h)\,dh\biggr)
  \prod_s \e_{m_{s}}(\rho_s)\,d\rho_r\dots d\rho_1.
  \label{eq:Z-r-k}
\end{multline}
We apply the definition (\ref{eq:eta}) of $\eta^k_{c,d}$
and move a finite summation 
$\sum_{\wt\Omega(r,k)}$ in the front of our integral.
Thus we get
$$
Z_{r,k}=\sum Y^\sigma_{r,k},
$$
where
\begin{multline}
Y^\sigma_{r,k}=
\sum_{c_1\in \Z} \dots \sum_{c_{n-2r}\in \Z}\,\,
\sum_{m_1\in \Z} \int\limits_{\rho_1\in\R} \dots
\sum_{m_r\in \Z} \int\limits_{\rho_r\in\R} \,\,\int\limits_{H_k}
\V_k(\partial)\Xi_k f(h)
\times
\\ \times
  \prod_{[c_{a_j}, \phi_\alpha]\in \sigma} e^{i c_{a_j}\phi_\alpha}
  \cdot 
  \prod_{[c_{b_l},\psi_\beta]\in\sigma} e^{i c_{b_l} \psi_\beta}
  \times\\ \times
  \prod_{[c_{a_j}, z_\gamma, c_{b_l}]\in \sigma}
  \Bigl(\sgn (t_\gamma) \exp\{-|c_{a_j}-c_{b_l}|\cdot |t_j|+ i(c_{a_j}+c_{b_l})\theta_j \}  \Bigr)
  \times\\ \times
  \prod_{[[d_s,z_\gamma]]\in\sigma}  e^{im_s\theta_\gamma} \bigl(e^{i\rho_s \theta_\gamma}- 
  e^{-i\rho_s \theta_\gamma}\bigr)
   \prod_s \e_{m_{s}}(\rho_s)\,dh\,d\rho_r\dots d\rho_1.
   \label{eq:Y-sigma}
\end{multline}

\sm

{\bf\punct Summation of distributions. Preliminaries.}
Lemmas below  follow Hirai \cite{Hir},
 but we need some additional details.

\begin{llemma}
	\label{l:cth}
	Let $f(t)$ be a smooth compactly supported
	  function
	on $\R$. Then
\begin{align}
\int_\R \int_\R f(t) e^{i\rho t} \,dt \coth(\pi\rho/2)\,d\rho
=\int_{\R} f(t) 2i \coth(t)\,dt;
\label{eq:cth}
\\
\int_\R \int_\R f(t) e^{i\rho t} \,dt \tanh(\pi\rho/2)\,d\rho= 
\int_{\R} f(t)\cdot
\frac{2i\,dt}{\sh t}.
\label{eq:th}
\end{align}
\end{llemma}

We will apply this lemma for odd functions $f(t)$. In this case,
the integrals in the right hand sides and
the repeated integrals in the left hand sides 
are absolutely convergent.

\sm

{\sc Proof.} We must evaluate Fourier transforms of tempered distributions
$\tanh(\pi\rho/2)$, $\coth (\pi\rho/2)$.
See Tables of Fourier transforms of distributions in \cite{BP}, Table 1, lines (396), (397), in the second case
we also must  apply
\cite{Bat}, formula (1.7.11).
Also, it is possible to apply 
formulas \cite{Bat-I}, (2.9.7), (2.9.8)
and continuity of the Fourier transform 
in the space of tempered distributions.
\hfill $\square$


 
 \begin{llemma}
 	\label{l:th-delta}
 	Let $f(t,\theta)$ be a smooth function with compact support,
 	satisfying $f(-t,\theta)=f(t,\theta)$. Then
\begin{multline} 
  \frac 12 \sum_{m\in \Z}\int_0^\infty\int_{-\pi}^{\pi}\int_{-\infty}^\infty f(t,\theta) e^{im\theta}
  \bigl(e^{i\rho t}-e^{-i\rho t} \bigr)\,dt\,d\theta\,\,
  \e_m(\rho)\,d\rho 
    =\\=
  \sum_{m\in \Z}\int_0^\infty\int_{-\pi}^{\pi}\int_{-\infty}^\infty f(t,\theta) e^{im\theta}
  e^{i\rho t}\,dt\,d\theta\,\,
  \e_m(\rho)\,d\rho
=\\=
  \Bigl\la \!\!\Bigl\la  f(t,\theta), 
   \pi\bigl( \coth(t/2)  \delta(\theta)+ \tanh (t/2)\delta(\theta-\pi)\bigr)
   \Bigr\ra \!\!\Bigr\ra _{(\R/2\pi \Z\,\times \R_+)}
   .
 \end{multline}
 \end{llemma}
 
 {\sc Proof.}
 The first equality follows from $f(-t,\theta)=f(t,\theta)$.
  The convergence of the series is obvious, and we 
 fulfill a formal calculation with distributions.
We have
$$
\sum_{m\in\Z} e^{i m\theta} =2\pi \delta(\theta)
$$
 and therefore
\begin{equation}
\sum_{l\in\Z} e^{i 2l\theta}=\pi \bigl(\delta(\theta)+\delta(\theta-\pi)\bigr), \qquad
\sum_{l\in\Z} e^{i (2l+1)\theta}=\pi \bigl(\delta(\theta)-\delta(\theta-\pi)\bigr).
\label{eq:delta-delta}
\end{equation}
 Keeping in the mind Lemma \ref{l:cth}, we evaluate 
 \begin{multline*}
 \sum_{m\in \Z}\int_{\rho\in\R} e^{im\theta}
 e^{i\rho t}
 \e_m(\rho)\,d\rho=\\=
 -\tfrac i2 \sum_{l\in \Z} e^{2il \theta}\cdot  \int\limits_{\rho\in\R} e^{i\rho t} \coth(\pi\rho/2)\,d\rho
 - \tfrac i2 \sum_{l\in \Z} e^{i(2l+1) \theta}\cdot \int\limits_{\rho\in\R} e^{i\rho t} \tanh(\pi\rho/2)\,d\rho
 =\\=
 -\tfrac {\pi i}2 \bigl(\delta(\theta)+\delta(\theta-\pi)\bigr)\cdot 2i \coth (t)
 -\tfrac {\pi i}2 \bigl(\delta(\theta)-\delta(\theta-\pi)\bigr)\cdot  \tfrac {2i}{ \sh(t)}
 =\\=
 \pi\bigl( \coth(t/2)  \delta(\theta)+ \tanh (t/2)\delta(\theta-\pi)\bigr).
  \qquad\quad \square
 \end{multline*}
 
 \sm

\begin{llemma}
	\label{l:tht-delta-2}
	Let $f(t,\theta)$ satisfy the conditions of the previous
	lemma and be supported by the set $|t|\le R$. Then
	
	\sm
	
	{\rm a)} For any $N>0$ there exists a constant $C(f,N)$
	such that for any $a$, $b\in \Z$,
	\begin{multline}
	\biggl|	\int_{0}^\infty f(t,\theta)
		 \exp\bigl\{-|a-b|\cdot |t|+ i(a+b)\theta\bigr\}\,dt\,d\theta
		\biggr|\le
		\\\le
		 \frac{C(f,N)}{(1+(a-b)^2) (1+|a+b|^N) } ,
		\label{eq:estimate1}
	\end{multline}
	where $C(f,N)$ admits a uniform estimate
	in terms of $R$ and numbers
	$$
	\max \biggl|\frac{\partial^{j} f}{\partial \theta^j}\biggr|,\qquad
	\max \biggl|\frac{\partial^{j+1} f}{\partial t\, \partial \theta^j}\biggr|,
	\qquad \text{where $j=0,\dots, N$.}
	$$
	
	{\rm b)} The following identity holds:
		\begin{multline}
	\sum_{a,b\in\Z}
	\int_{-\pi}^{\pi}
	\int_{-\infty}^\infty f(t,\theta)
	\sgn (t) \exp\bigl\{-|a-b|\cdot |t|+ i(a+b)\theta\bigr\}\,dt\,d\theta
	=\\=
	\Bigl\la\!\!\Bigl\la f(t,\theta)
	\pi
	\bigl(\coth(t/2)\delta(\theta)+ \tanh(t/2)\delta(\theta-\pi)\bigr)
	\Bigr\ra\!\!\Bigr\ra_{(\R/2\pi\Z\times \R)}
	,
	\label{eq:sum-dist}
	\end{multline}
\end{llemma}

{\sc Proof.} a) We integrate by parts $N$ times with respect to
$\theta$ and one time with respect  to $t$.

\sm

b) By a) the double series in the left-hand side of
(\ref{eq:sum-dist}) absolutely converges. 
We change summation indices to
$m=a$, $n=a-b$. Thus we must evaluate a sum of distributions
$$
\sum_{n\in \Z} \sum_{m\in \Z} \exp\bigl(-|n| t\bigr) 
\exp\bigl(i(2m-n)\theta\bigr),
$$
to simplify notation we assume $t>0$, recall that we consider functional on the space  functions that odd in $t$.
By (\ref{eq:delta-delta}) we get
\begin{multline*}
\sum_{k\in \Z} \exp\bigl(-|2k|\,t\bigr)
\cdot \pi \bigl(\delta(\theta)+\delta(\theta-\pi)\bigr)
+
\sum_{k\in \Z} \exp\bigl(-|2k+1|\,t\bigr)
\cdot
\pi \bigl(\delta(\theta)-\delta(\theta-\pi)\bigr)
=\\=
\coth t \cdot \pi \bigl(\delta(\theta)+\delta(\theta-\pi)\bigr)
+
\frac 1{\sh t}\cdot \pi \bigl(\delta(\theta)-\delta(\theta-\pi)\bigr)
\end{multline*}
The last pass is a formal manipulation with series,
it is justified by the dominated convergence 
of the expression
$$
\sum_{k\in \Z} \int_0^\infty f(t,0) e^{-|n|\,t}\,dt
= \int_0^\infty f(t,0) \Bigl(1+ \frac{2e^{-t}}{1-e^{-t}}
\Bigr) \,dt,
$$
recall that $f(t,\theta)=0$.
After a simple transformation we come to
$$
\qquad\qquad\qquad
\pi \coth (t/2)\, \delta(\theta)+\pi \tanh (t/2)\, \delta(\theta-\pi).
\qquad\qquad\qquad\square
$$

\begin{llemma}
	\label{l:delta}
Let $f(\phi)$ be a continuous piece-wise smooth function
on $\R/2\pi\Z$.	Then there exist a constant $C(f)$ 
such that
\begin{equation}
\biggl|\int_{-\pi}^{\pi} f(\phi) e^{-ik\phi}\,d\phi\biggr|
\le\frac {C(f)} {(1+|n|^2)},
\label{eq:estimate2}
\end{equation}	
where $C(f)$ can be estimated in terms of 
$$
\max|f|,\qquad \max|f'|,\qquad \max|f''|
$$
and the number of singular points of $f$. Moreover
$$
\sum_{k\in \Z} \int_{-\pi}^{\pi} f(\phi) e^{-ik\phi}\,d\phi
= 2\pi\, f(0).
$$
\end{llemma}

{\sc Proof.} To obtain (\ref{eq:estimate2}) we integrate
two times by parts (after the first integration boundary
term do not appear). This implies the absolute convergence of the Fourier series and the second statement.
\hfill $\square$


\begin{llemma}
	\label{l:convergence-1}
	Let $f(\phi,t,\theta)$ be a continuous function
	piece-wise smooth function on
	 $\R/2\pi \Z \times \R/2\pi\Z$, which is smooth for 
	 any fixed $\phi$.
	 Then the following double sum
	 absolutely converges
	 \begin{equation}
	 	\sum_{a\in\Z} \sum_{c\in \Z}
	 	\int_{-\pi}^{\pi}\int_{-\pi}^{\pi}
	 	f(\phi,t,\theta)
	 	\sgn (t) \exp\bigl\{-|a-b|\cdot |t|+ i(a+b)\theta\bigr\}\,d\theta\,d\phi
	 \end{equation}	 
\end{llemma}

{\sc Proof.} This follows from estimates
(\ref{eq:estimate1}), (\ref{eq:estimate2}).
\hfill $\square$

\begin{llemma}
	\label{l:convergence-2}
Let $f(t_1,t_2 ,\theta_1,\theta_2)$
be a smooth function on $\R\times\R\times\R/2\pi\Z\times\R/2\pi\Z$
 odd in
$t_1$ and odd in $t_2$. Then the following
series absolutely converges
\begin{multline*}
	\sum_{a\in\Z}  \sum_{c\in \Z}
	\int_{-\pi}^{\pi}\int_{-\pi}^{\pi}
	f(t_1,t_2 ,\theta_1,\theta_2)
	\sgn (t) \exp\bigl\{-|a-b|\cdot |t_1|+ i(a+b)\theta_1\bigr\}
		\times\\\times
	\sgn (t) \exp\bigl\{-|c-d|\cdot |t_2|+ i(c+d)\theta_2\bigr\}	
	\,d\theta_1\,d\theta_2.
\end{multline*}
\end{llemma}	

{\sc Proof.} This follows from (\ref{eq:estimate1}).
\hfill $\square$

\sm

{\bf\punct Summation of distributions.%
	\label{ss:summation}}
%
Formally transforming (\ref{eq:Y-sigma})
we get 
\begin{multline}
 Y^\sigma_{r,k}=
 \biggl\la \!\!\!\biggl\la
 \Updelta_k(\partial)\Xi_k f(\phi,\psi,t,\theta)\,,\,
 \prod_{\alpha=1}^p \bigl(\sum_{a\in \Z} 
 e^{ia\psi_\alpha}\bigr)
 \cdot \prod_{\beta=1}^q \bigl(\sum_{b\in \Z} 
 e^{ib\phi_\beta}\bigr)
 \times\\\times
 \prod_{\gamma:\, [c_{a_j}, z_\gamma,c_{b_l}]\in\sigma}
 \Bigl(
 \sum_{a,b\in\Z} \sgn (t_\gamma) 
 \exp\bigl\{-|a-b|\cdot |t_\gamma|+ i(a+b)\theta_\gamma\bigr\}
 \Bigr)
 \times\\\times
 \prod_{\gamma:\, [[d_s, z_\gamma,]]\in\sigma}
  \Bigl( \sum_{m\in \Z} e^{im\theta_\gamma}
 \int_{\rho\in \R}
 e^{i\rho t_\gamma}\bigl(e^{i\rho t_\gamma}-e^{-i\rho t_\gamma} \bigr) 
 \e_m(\rho)\,d\rho\Bigr)
 \biggr\ra \!\!\!\biggr\ra_{(H_k)}
 \end{multline}

Applying Lemmas \ref{l:th-delta}--\ref{l:delta} we come to
\begin{multline}
Y^\sigma_{r,k}=
2^{n-2k+r}\pi^{n-k}
\times\\\times
\prod_{\alpha=1}^p \delta(\phi_\alpha) \cdot \prod_{\beta=1}^q \delta(\psi_\beta)
\cdot \prod_{\gamma=1}^k\bigl( \coth(t_\gamma)  \delta(\theta_\gamma)+ \tanh (t_j)\delta(\theta_\gamma-\pi)\bigr).
\label{eq:Y-sigma-2}
\end{multline}
and we observe that for all $\sigma$ the sums  $Y^\sigma_{r,k}$ are equal.
However, in the initial expression for
$Y^\sigma_{r,k}$ we have a successive integration-summation,
the calculation above assumes  changings of the order
of the summation.

  \sm
  
 {\sc Step 1.} We can change a summation in $m_r$
  with any integration in $d\phi_\alpha$, $d\psi_\beta$,
  and $dt_\gamma$, $d\theta_\gamma$ except $\gamma$ linked with $r$
  on the diagram $\sigma$.
    For this, we use two identities, the first is
    \begin{equation*}
    \sum_{m\in \Z} \Bigl(\int_{-\pi}^\pi \int_{-\pi}^\pi e^{im\theta} F(\theta, \phi)\,d\theta\,d\phi\Bigr) 
    =
    \int_{-\pi}^\pi \Bigl( \sum_{m\in \Z} \int_{-\pi}^\pi e^{im\theta} F(\theta, \phi)\,d\theta\Bigr) \,d\phi
    \end{equation*}
for a  function $F$, which is piece-wise 
smooth on a torus.
The second is 
$$
    \sum_{m\in \Z} \Bigl(\int_{-\pi}^\pi \int_\R e^{im\theta} F(\theta, \tau)\,d\theta\,d\tau\Bigr) 
    =
    \int_\R \Bigl( \sum_{m\in \Z} \int_{-\pi}^\pi e^{im\theta} F(\theta, \tau)\,d\theta\Bigr) \,d\tau
    $$
    for a smooth compactly supported $G$ on $\R/2\pi\Z\times \R$.

  The same property holds for the integration in $d\rho_r$,

  So we can start  the integration-summation
  from the integral
  \begin{multline*}
  \sum_{m_r} \int_{\rho_r\in\R} \int_{t\in\R}
  \int_{\theta\in[0,2\pi]} \V_k(\partial) \Xi_k f(\phi,\psi,t,\theta)
  \times\\ \times
   e^{im_r\theta_\gamma}
  \bigl(e^{i\rho_r t_\gamma}-e^{-i\rho_r t_\gamma}
   \bigr)\,dt_\gamma\,d\theta_\gamma\,\,
  \e_m(\rho_r)
  \,d\rho_r
  \end{multline*}
  We apply Lemma \ref{l:th-delta} and continue the process.
  In this way, we  perform successive integrations and summations
  in (\ref{eq:Y-sigma}) one after another for all $\sum_{m_s}\int d\rho_s$.
  
  \sm
  
   {\sc Step 2.} We start a successive summation with respect to
  $c_{n-2r}$, $c_{n-2r-1}$, \dots. As above, we can change one summation in $c_j$ with an integration
  in $d\phi_\alpha$, $d\psi_\beta$, $dt_\gamma$, $d\theta_\gamma$
  (if this variables are not linked with $c_a$ in the diagram $\sigma$.)
If we meet $c_a$ linked with $\phi_\alpha$ or $\psi_\alpha$,
we apply Lemma \ref{l:delta}.
For the factors of the type
$$
\Bigl(\sgn (t_\gamma) \exp\{-|c_{a_j}-c_{b_l}|\cdot |t_j|+ i(c_{a_j}+c_{b_l})\theta_j \}  \Bigr)
$$
the corresponding summations
$\sum_{c_{a_j}}$, $\sum_{c_{b_j}}$
generally are not adjacent in (\ref{eq:Y-sigma}) and we can not immediately apply Lemma \ref{l:th-delta}. But Lemmas \ref{l:convergence-1},
\ref{l:convergence-2} allow to change adjacent summations
$\sum_{c_j}$, $\sum_{c_{j+1}}$ in two case

\sm

--- if both $c_j$, $c_{j+1}$ are connected with black boxes;

\sm

--- if precisely one of  $c_j$, $c_{j+1}$ is connected with
a black box.

\sm

This allows a consequent application of Lemmas 
\ref{l:delta} and \ref{l:tht-delta-2}.

\sm

In this way, we justify (\ref{eq:Y-sigma-2} and get
$$
Z_{r,k}=\# \wt \Omega(r,k)\cdot Y^\sigma_{r,k}
$$
This implies Theorem 1.

\sm

\section{Refinements of the orthogonal decomposition}

\COUNTERS

\begin{figure}
$$
\epsfbox{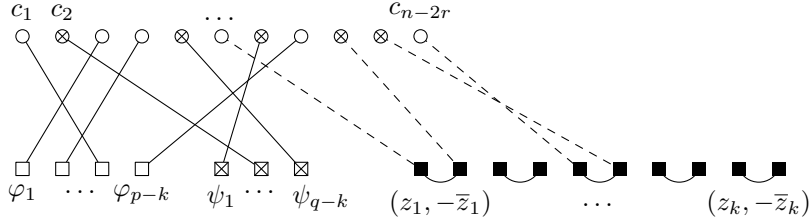}
$$
\caption{Elements of $\Omega^\circ(A|r,k)$.\label{fig:last}}
\end{figure}

{\bf \punct The decomposition with respect to the parameters $\mathbf{(A;c)}$.%
	\label{ss:th2}}
Fix $r$ and $(A;c)$. 
Denote  by $L_{A;c}$ the subspace in the Plancherel decomposition (\ref{eq:decomposition-u}),
which is the integral of all representations with given $(A;c)$. To write the projector to $L_{A;c}$ we need some  notations.
Denote by $\Omega^\circ(A|r,k)$ the set of all diagrams of the form 
shown on Figure \ref{fig:last}. These diagrams are obtained from elements of $\Omega(A|r,k)$
by forgetting black circles and adjacent arcs. These diagrams split into elements of the following 4 types

a) Arcs $\zcircle$---$\zsquare$ or  $c_{a_j}$---$\phi_\alpha$, where
$a_j\in A$.

\sm

b) Arcs $\zotimes$---$\zboxtimes$ or $c_{b_l}$---$\psi_\beta$, 
where $b_l\in B$.

\sm

c) Chains $\zcircle$---$\zblacksquare$---$\zblacksquare$---$\zotimes$ 
or $c_{a_j}$---$z_\gamma$---$(-\ov z_\gamma)$---$c_{b_l}$.

\sm

d)$^*$ Arcs $\zblacksquare$---$\zblacksquare$, or $z_\gamma$---$(-\ov z_\gamma)$.

\begin{ttheorem} 
 An invariant distribution $\Theta_{A;c}$ determining an orthogonal  projector
 to $L_{A;c}$ is given by
 the formula
 \begin{multline*}
M_*  \lla f, \Theta_{A;c}\rra_{(G)}= \\=
(-1)^{n(n-1)/2+pq+qr+r(r-1)/2} 
\sum_{k=r}^n \frac{ 2^{n-2k}\pi^{n-k} \omega_k} {r!}
\lla \V_k(\partial) \Xi_k f, \zeta^k_{A;c}\rra_{(H_k)}, 
 \end{multline*}
 where
 \begin{multline*}
 \zeta^k_{A;c}= 
 \sum_{\sigma\in \Omega^\circ(A|r,k)} 
 \prod_{[c_{a_j}, \phi_\alpha]\in \sigma} e^{i c_{a_j}\phi_\alpha}
 \cdot
 \prod_{[c_{b_l},\psi_\beta]\in\sigma} e^{i c_{b_l} \psi_\beta}
 \times\\ \times
  \prod_{[c_{a_j}, z_\gamma, c_{b_l}]\in \sigma}
\Bigl(\sgn (t_\gamma) \exp\{-|c_{a_j}-c_{b_l}|\cdot |t_j|+ i(c_{a_j}+c_{b_l})\theta_j \}  \Bigr)
  \times  \\\times 
 \prod_{[[z_\gamma, -z_\gamma]]\in\sigma} 
 \bigl( \coth(t_\gamma/2)  \delta(\theta_\gamma)+ \tanh (t_\gamma/2)\delta(\theta_\gamma-\pi)\bigr).
 \end{multline*}
\end{ttheorem}

 A calculation of the distributions $\zeta^k_{A;c}$
 are the same as above. 
 It is important that Theorem 2 from \cite{Hir}
 allows the integration by parts for functions $\eta^k_{A;c,d}$ defined in
 Subsect. \ref{ss:transformations}.
 We symmetrize with respect to $S_r\ltimes \Z_2^r$ instead of
 $W^{(r)}$ in (\ref{eq:Z-r-k}). In the calculation described 
 in Subsect. \ref{ss:summation}  we make only Step 1.

 \sm
 
 {\bf \punct The decomposition with respect to the parameters  $\mathbf{(A;c)}$ and $\mathbf{m}$.%
 	\label{ss:th3}}
 Take $r$,  $A$, $(c_1, \dots ,  c_{n-2r})$, and
 \begin{equation}
 ( m_1+i\rho_1,\,\dots,m_r+i\rho_r)
 \label{eq:m+irho}
 \end{equation}
 Consider 
 a representation $T_{A;c,m,\rho}$ of $\U(p,q)$ unitary induced from a representation
 (\ref{eq:rep-parabolic}).
 The hyperoctahedral group $S_r\ltimes \Z_2^n$ acts on the set of collections 
 (\ref{eq:m+irho}) by permutations and complex cojugations $m_s+i\rho_s\mapsto m_s-i\rho_s$.
Elements of this group send representations $T_{A;c,m,\rho}$  to equivalent representations.
In particular, we can assume that
\begin{equation}
m_1\ge m_2\ge\dots\ge m_r.
\label{eq:mmm}
\end{equation}
Denote by $L_{A;c,m}$ the integral of all representations with fixed
$A$, $c$, $m$ in $L^2(G)$. We intend to write a projector to
 $L_{A;c,m}$. 

Sinse a collection (\ref{eq:mmm})  can contain repeating entries, we 
 will use an alternative notation for the same  collection,
\begin{equation}
\underbrace{\ov{\ov m}^{\,1}, \dots, \ov{\ov m}^{\,1}}_{\text {$u_1$ times}} >
\underbrace{\ov{\ov m}^{\,2}, \dots, \ov{\ov m}^{\,2}}_{\text {$u_2$ times}}>
\underbrace{\ov{\ov m}^{\,3}, \dots, \ov{\ov m}^{\,3}}_{\text {$u_3$ times}}>\dots
\label{eq:m-new}
\end{equation}
(here $u_\iota>0$, $\sum u_\iota=r$).

Define function $\frf_m(t)$, where $m\in\Z$, by
$$
\frf_m(t)=
\begin{cases}
 \coth(t) ,\qquad \text{where $m\in 2\Z $;}
 \\
1/\sh (t),\qquad \text{where $m\in 2\Z+1 $.}
\end{cases}
$$

\begin{ttheorem}
 The invariant distribution $\Theta_{A;c,m}$ determining the orthogonal projector
 to $L_{A;c,m}$ is given by 
 \begin{equation*}
 M_* \lla f,\Theta_{A;c,m}\rra_{(G)}
  = \frac{(-1)^{n(n-1)/2+pq+qr+r(r-1)/2}} {2^r} \sum_{k=r}^n \lla\V_k(\partial) \Xi_k f, \zeta_{A;c,m}\rra_{(H_k)},
 \end{equation*}
 where
 \begin{multline*}
 \zeta_{A;c,m}=
 \frac{1}{\prod u_\iota!} \sum_{\sigma\in \Omega(r,k)}
      \prod_{[c_{a_j}, \phi_\alpha]\in \sigma} e^{i c_{a_j}\phi_\alpha}
 \cdot 
 \prod_{[c_{b_l},\psi_\beta]\in\sigma} e^{i c_{b_l} \psi_\beta}
 \times\\ \times
 \prod_{[c_{a_j}, z_\gamma, c_{b_l}]\in \sigma}
\Bigl(\sgn (t_\gamma) \exp\{-|c_{a_j}-c_{b_l}|\cdot |t_j|+ i(c_{a_j}+c_{b_l})\theta_j \}  \Bigr)
 \times\\ \times
 \prod_{[[d_s,z_\gamma]]\in\sigma}  e^{im_s\theta_\gamma} \frf_{m_s}(t_\gamma).
 \end{multline*}
\end{ttheorem}

{\sc Proof.} In the Plancherel formula we have summation over the set 
$m_1$, \dots, $m_r\in \Z$, $\rho_1>\rho_2>\dots>\rho_r$. We can replace this domain
by any fundamental domain of the hyperoctahedral group $S_r\ltimes \Z_2^r$.
Denote the parameters $\rho$ corresponding to (\ref{eq:m-new})
by
$$
\rho_1^1,\rho_1^2,\dots,\rho_1^{u_1}, \rho_2^1,\rho_2^2,\dots,\rho_2^{u_2},\rho_3^1,\rho_3^2,\dots,\rho_3^{u_3},\dots
$$
We choose a fundamental domain determined by
(\ref{eq:m-new}) and
$$
\rho_1^1>\rho_1^2>\dots>\rho_1^{u_1}>0,\qquad
\rho_2^1>\rho_2^2>\dots>\rho_2^{u_2}>0,\qquad 
\dots
$$
Now problem is reduced to an evaluation of
$$
\int\limits_{\begin{matrix}
       \rho_1^1>\rho_1^2>\dots>\rho_1^{u_1}>0,\\
       \rho_2^1>\rho_2^2>\dots>\rho_2^{u_2}>0,\\
       \hdotsfor{1}
      \end{matrix}
} \eta_{A;c,d} \prod_{s}\e_m(\rho_s)\,d\rho.
$$
Using symmetry, we change this to
$$
 \frac{1}{2^r\prod u_\iota!}
\int\limits_{\begin{matrix}
       (\rho_1^1,\rho_1^2,\dots,\rho_1^{u_1})\in\R^{u_1},\\
       (\rho_2^1,\rho_2^2,\dots,\rho_2^{u_2}) \in\R^{u_2},\\
       \hdotsfor{1}
      \end{matrix}
} \eta_{A;c,d} \prod_{s}\e_m(\rho_s)\,d\rho.
$$
We pass to the sum $\sum_{\sigma\in\Omega(A|r,k)}$
and integrate it termwise in $\rho_r$, \dots $\rho_1$ using 
(\ref{eq:th}), (\ref{eq:cth}).
\hfill $\square$

\tt
\noindent
Math. Dept., University of Vienna; \\
Institute for Theoretical and Experimental Physics (Moscow) \\
MechMath Dept., Moscow State University\\
Institute for Information Transmission Problems\\
e-mail:neretin(frog)mccme.ru\\
URL: http://mat.univie.ac.at/~neretin/


\begin{thebibliography}{cc}
	
	\bibitem{AU}
	Alldridge, A.; Upmeier, H.
	{\it Toeplitz operators on Hardy spaces over $\SL(2,\R)$: irreducibility and representations.} In
	{\it Geometry and analysis on finite- and infinite-dimensional Lie groups}, 173--209,
	Banach Center Publ., 55, Polish Acad. Sci., Warsaw, 2002. Barut, Asim O.
	
	
	\bibitem{BR}
	Barut A. O; Raczka, R.
{\it	Theory of group representations and applications.}
PWN,	Warszawa, 1977. 
	
	
	
	\bibitem{BS}
	Ben Said, S.
{\it	Espaces de Bergman pond\'er\'es et s\'erie discr\'ete holomorphe de $\wt{\U(p,q)}$.}
	J. Funct. Anal.,
	173,   (2000), 1,  154-181.
	

	\bibitem{BP}
	Brychkov, Yu. A.; Prudnikov, A. P. {\it Integral transforms of generalized functions.}
	Nauka, Moscow, 1977 (Russian); English transl. from the Second Russian edition:
	Gordon and Breach, New York, 1989.
	
	\bibitem{Bat}
Erd\'elyi, A.; Magnus, W.; Oberhettinger, F.; Tricomi, F. G. {\it Higher transcendental functions.
Vols. I. Based, in part, on notes left by Harry Bateman.}
McGraw-Hill Book Company, Inc., New York-Toronto-London, 1953.	

\bibitem{Bat-I}
Erd\'elyi, A.; Magnus, W.; Oberhettinger, F.; Tricomi, F. G. {\it Tables of integral transforms. Vol. I.
Based, in part, on notes left by Harry Bateman.} McGraw-Hill Book Company, Inc., 
New York-Toronto-London, 1954. 
	
	\bibitem{FL}
	Frenkel I., Libine L.
	{\it Split quaternionic analysis and separation of the series for $\SL(2,\R)$ and $\SL(2,\C)/\SL(2,\R)$}.
	Adv.  Math.
	228 (2011), 2, 678-763.
	
	\bibitem{GG}
	Gelfand, I. M.; Gindikin, S. G. {\it Complex manifolds whose spanning trees are real semisimple Lie groups,
	and analytic discrete series of representations.}
	Funct. Anal.  Appl., 1977, 11:4, 258-265.
	
	\bibitem{GN} Gelfand I.M., Naimark M.A.
	{\it Unitary representations of the classical
		groups.} Trudy Mat. Inst. Steklov., vol. 36, 1950; German transl.: Gelfand,
	I. M., Neumark, M. A. {\it Unit\"are Darstellungen der klassischen Gruppen.}
	Akademie-Verlag, Berlin, 1957.
	
	
	
	\bibitem{Gin1}
	Gindikin S., {\it Conformal analysis on hyperboloids,} J. Geom. Phys., 10, 175-184 (1993). 
	
	\bibitem{Gin2}
	Gindikin, S.
	{\it Integral geometry on $\SL(2;\R)$}. 
	Math. Res. Lett. 7 (2000), no. 4, 417-432. 
	
	\bibitem{Gin3}
	Gindikin, S. {\it  An analytic separation of series of representations for $\SL(2;\R)$}.
	Mosc. Math. J. 2 (2002), no. 4, 635-645.
	
	\bibitem{GK}
	Gindikin, S, Kr\"otz, B., \'Olafsson, G.
	{\it Hardy spaces for non-compactly causal symmetric spaces and the most continuous spectrum.}
	Math. Ann. 327 (2003), no. 1, 25-66.
	
	
	
	
	\bibitem{Harish}
	Harish-Chandra,
	{\it A formula for semisimple Lie groups.}
	Amer. J. Math. 79 (1957) 733-760.
	
	\bibitem{Harish-discrete}
	Harish-Chandra,
{\it Discrete series for semisimple Lie groups. II. Explicit determination of the characters.}
Acta Math. 116 (1966) 1-111. 

\bibitem{Harish-Plancherel}
Harish-Chandra
{\it Harmonic analysis on real reductive groups. III. The Maass-Selberg relations and the Plancherel formula.}
Ann. of Math. (2) 104 (1976), no. 1, 117-201. 
	
\bibitem{HW}
Herb, R., Wolf, J. A.,
{\it The Plancherel theorem for general semisimple groups.}
Compositio Math. 57 (1986), no. 3, 271-355. 
	
	
	\bibitem{HO}
	 Hilgert J.,  \'Olafsson G., {\it Causal symmetric spaces. geometry and harmonic analysis,} 
	Academic Press Inc., San Diego, CA, 1997.
	
	\bibitem{Hir}
	Hirai, T.
{\it The Plancherel formula for $\SU(p,q)$.} 
J. Math. Soc. Japan 22 (1970), 134-179. 
	
	
	\bibitem{Kna}
	Knapp A.
{\it	Representation theory of semisimple groups : An overview based on examples},  Princeton: University Press, 2001.
	
	\bibitem{KM}
	Kobayashi T., Gen Mano
	{\it The inversion formula and holomorphic extension of the minimal representation of the conformal group,}
	in {\it Harmonic Analysis, Group Representations, Automorphic Forms and Invariant Theory:
	In Honour of Roger E. Howe,} World Scientific, 2007, pp. 159-223.
	
	\bibitem{Kof}
	Koufany, K.; \O rsted, B. {\it Function spaces on the Olshanski semigroup and the Gelfand-Gindikin program.}
	Ann. Inst. Fourier (Grenoble) 46 (1996), no. 3, 689-722. 
	
	\bibitem{Mol1}
	Molchanov, V. F.
	{\it Quantization on the imaginary Lobachevski plane}. 
	Funct. Anal. and Appl., 1980, 14:2, 142-144
	
	
	\bibitem{Mol2}
	Molchanov, V. F.
	{\it Separation of series for hyperboloids.}
	Funct. Anal. Appl. 31 (1997), no. 3, 176-182. 
	
	\bibitem{Ner-compl}
	Neretin, Yu. A.
	{\it Discrete occurrences of representations of the complementary series 
	in tensor products of unitary representations.}  Functional Anal. Appl. 20 (1986), no. 1, 68-70. 
	
	\bibitem{Ner-rest}
	Neretin, Yu. A.
	{\it Restriction of functions holomorphic in a domain to curves lying on the boundary of the domain, 
	and discrete $\SL_2(\R)$-spectra}. 
	Izv. Math. 62 (1998), no. 3, 493-513.
	
	\bibitem{Ner-separ}
	Neretin, Yu. A.
	{\it On the separation of spectra in the analysis of Berezin kernels}. 
	Funct. Anal. Appl. 34 (2000), no. 3, 197-207 
	
	\bibitem{NO}
	Neretin, Yu. A.; Olshansky, G. I. {\it  Boundary values of holomorphic functions, 
	singular unitary representations of the groups $\OO(p,q)$ and their limits as 
		$q\to\infty$.} (Russian) Zap. Nauchn. Sem. POMI, 223 (1995), 9-91;
	translation in J. Math. Sci. (New York) 87 (1997), no. 6, 3983-4035 
	
	\bibitem{Ols1}
	Olshanski, G. I. {\it Complex Lie semigroups, Hardy spaces and the Gelfand-Gindikin program.}
	(Russian) Problems in group theory and homological algebra, 85-98, 161, Yaroslav. Gos. Univ., Yaroslavl, 1982; English transl.:
	Differential Geom. Appl. 1 (1991), no. 3, 235-246.
	
	
	
	
	
	
	\bibitem{War}
	Warner, G. {\it Harmonic analysis on semi-simple Lie groups. II.}
	Springer-Verlag, New York-Heidelberg, 1972.
	
\end{thebibliography}
\end{document}